# cK¢, A MODEL TO REASON ON LEARNERS' CONCEPTIONS

Nicolas Balacheff
CNRS - Laboratoire d'Informatique de Grenoble
Nicolas.balacheff@imag.fr

*Understanding learners' understanding is a key requirement for an efficient design of teaching situations and learning environments, be they digital or not. This keynote outlines the modeling framework cK¢ (conception, knowing, concept) created with the objective to respond to this requirement, with the additional ambition to build a bridge between research in mathematics education and research in educational technology. After an introduction of the rationale of cK¢, some illustrations are presented. Then follow comments on cK¢ and learning. The conclusion evokes key research issues raised by the use of this modeling framework.*

**A short story**

The model I will introduce has a long history. Its construction started at the end of the 80s with the project of bridging artificial intelligence and didactics of mathematics, and the objective of enhancing the design of computer-based learning environments. The design of these environments had a ternary structure including a model of the learner, a model of the content to be learned and an instructional model. These models are nowadays either effectively implemented, or only involved in the design phase. For the last two, research has constantly been very active with some promising progress. On the contrary, modeling the learner proved to be a real challenge, and actually it is still the case despite expectations raised by recent research on educational data mining and learning analytics.

The theoretical framework within which I was working, the Theory of Didactical Situations (TDS), is based on ideas that some may consider as precursors of "learning games" now celebrated by researchers in educational technology:

> "*Modeling a teaching situation consists of producing a game specific to the target knowledge among different subsystems: the educational system, the student system, the milieu, etc.*"
> (Brousseau 1986/1997 p.47).

Within this approach, the teacher is "*a player faced with a system, itself built up from a pair of systems: the student and, let us say for the moment, a 'milieu' that lacks any didactical intentions with regards to the student*" (ibid. p.40). A systemic approach to learner modeling makes possible expressing learning as adaptation and adaptation as construction. However, while research went quite far in modeling didactical situations, the progress on modeling the learner "*subsystem*" was rather limited. Indeed, there was a lot of research within different frameworks, with different concepts and a variety of terms, reporting observations of learners' behaviors. But the distance between the content of these reports and models that we could use to inform the design of learning environments, was quite significant.

As a researcher in mathematics educations, I had another motivation. I started my research on the learning of mathematical proof in the 70s, with an exploration of the structure of proofs produced by students, using graphs to represent them. This approach was intrinsically limited, so it was not too difficult to understand the critiques of Guy Brousseau and Gérard Vergnaud. The former drew my attention to the role of the situation; the latter pointed to the cognitive complexity of producing proofs. Then, I engaged in a completely different direction which resulted at the end of 80s in a first picture of what could be the genesis of mathematical proof. While this picture could help designing didactical situations, it didn't shed light on the underlying mental processes. In other words, it was a didactical study, not a psychological study. Still, there was a weak point: it left open a gap between proving and knowing. To close this gap I had to find a way to model the learner's ways of knowing. The Vergnaud model of concept and conceptual field offered a possible solution:

> "[…] *from a developmental point of view, a concept is altogether: a set of situations, a set of operational invariants (contained in schemes), and a set of linguistic and symbolic representations.*" (Quoted from Vergnaud 2009 p.94, but this characterization goes back to the early 80's)

Vergnaud introduced the notation: $C=(S, I, \mathscr{S})$ in which the components refer respectively to each of the three sets mentioned above. He emphasized that these components cannot be separated; they have to be considered all at the same time when studying learners' development. This characterization has direct connections with the TDS description of the relation between a learner and a milieu based on different forms of knowledge (Brousseau ibid. p.61):

> [1] The models for action governing decisions.
> [2] The formulation of the descriptions and models.
> [3] The forms of knowledge which allow the explicit "control" of the subject's interactions in relation to the validity of her statements.

Apart from the set of situations (S) which is implicitly shared, common elements are related to action (I) and formulation ($\mathscr{S}$). Then, one element is missing which corresponds to the terms "control" used in the TDS description. I mentioned this lack when discussing the Vergnaud concept of theorem-in-action. A theorem, and the same applies to a theorem-in-action, is both a tool and a statement: "if A then B" is a tool to obtain B if A is valid, it is also a statement which has a truth value. This duality of "*the operational form and the predicative form of knowledge*", as Vergnaud expresses it, facilitated keeping implicit the control dimension in the characterization he proposed. However, after Polya and a long tradition of research on metacognition, Schoenfeld (1985 pp. 97-143) has shown the crucial role of control in problem-solving. The suggestion I made consists in introducing explicitly this dimension of control in Vergnaud model. This is the origin of the quadruplet I describe below.

Before entering the main content of this keynote, I would like to address an issue which led me use the term "conception" and not the term "knowledge" as it is classical in educational technology and mathematics education as well. Most of our research is based implicitly or not on the hypothesis that learners act as rational subjects. But, one often is faced to rational thinking co-existing with knowledge which looks contradictory (from the observer's point of view). Bourdieu (1990) proposed a solution to this paradox: "*The calendar thus creates ex nihilo a whole host of relations […] between reference-points at different levels, which never being brought face to face in practice, are practically compatible even if they are logically contradictory*" (ibid. p. 83). The key elements are time on one hand, and on the other hand the diversity of situations. Time organizes the subjects' decisions sequentially in such a way that even contradictory, they are equally operational because appearing at different periods of their history: contradictory decisions can ignore each other. The diversity of the situations introduces an element of a different type. It is a possible explanation insofar as one recognizes that each decision is not of a general nature but that it is related to a specific sphere of practice (we would say, nowadays, that it is situated) within which it is acknowledged as efficient. Within a sphere of practice learners are coherent and successful; they are non-contradictory.

Contradictions (and failures) appear when learners are faced with situations foreign to their sphere of practice but in which they have nevertheless to produce a response (e.g. a question from an interviewer). They mobilize what they have available which worked elsewhere, but more often than not this ends in systematically making errors. The classical position in the 80s was to consider these errors as symptoms of misconceptions. This term used to come with expressions like "naive theory", "private concepts", "beliefs" or even "mathematics of the child". Such views missed the fact that "*a child may not be 'seeing' the same set of events as a teacher, researcher or expert. [...] many times a child's response is labeled erroneous too quickly and [...] if one were to imagine how the child was making sense of the situation, then one would find the errors to be reasoned and supportable* " (Confrey 1990 p.29). Agreeing with this position, I renounced using the term "misconception". Still, recognizing that learners may have different models-in-action to mobilize for (what we consider as) the same piece of knowledge, I needed a term, but one different from "knowledge" because of the issue raised by the observation of possible contradictions in learners behaviors. A possible term was "conception" largely used in science education to denote theory-in-action. Most often than not conception functioned as a tool in discourses but it was not taken as an object of study as such (Artigue 1991, p.266), although there was an acknowledged need (e.g. Vinner 1983, 1987) for a better grounded definition of conceptions, and for tools allowing analyzing their differences and commonalties.  In the following section, I propose a definition of "conception", and then describe a model revisiting the Vergnaud's triplet.

**Behavior, conception and knowing**

The only indicators one has to get an insight into learners' understanding are their behaviors and products which are consequences of the conceptions they have engaged. Such evaluations are possible and their results are significant only in the case where one is able to establish a valid relationship between the observed behaviors and the invoked conception. This relation has been relatively "hidden" as such for a long while as a result of the fight against behaviorism, but it has always been present in educational research at least at the methodological level. Indeed, the key issue is that *the meaning of a piece of knowledge cannot be reduced to behaviors, whereas meaning cannot be characterized, diagnosed or taught without linking it to behaviors*.

Being a tangible manifestation of the relationships between a person and her environment, a "behavior" depends on the characteristics of this person as well as on the characteristics of her environment. A now well documented example is that of an instrument which at the same time facilitates action if the user holds the required competence, and on the other hand limits this action because of its own constraints (Rabardel 1995, Resnick & Collins 1994, p.7).

The words "person" and "environment", which I am using here, refer to complex realities whose aspects are not all relevant for our investigations. One may want to ignore the clothes the person wears and the shape of the room in which he or she stands (although we have always to be prepared to consider seriously features initially downplayed). What is of interest is the person from the point of view of his or her relationship to a piece of knowledge. For this reason I will refer from now on to the learner as a reduction, if I dare saying so, of the person to her cognitive dimension. In the same way, I do not consider the environment in all its complexity, but only those of its features that are relevant with respect to a given piece of knowledge. Actually, this corresponds to the TDS concept of *milieu*, which is a kind of projection of the environment onto its epistemic dimension: the milieu is the learner's antagonist system in the learning process (Brousseau, 1997 p.57)

This situatedness nature of a conception suggests not considering it as a property which can be ascribed only to the learner but as a property of the interacting system formed by the learner and his or her antagonist milieu, to which I will refer as the "learner/milieu system". What is requested for this property to be valid is that the system satisfies the necessary conditions for its viability. I mean that the system has the capacity to recover equilibrium after a perturbation which otherwise would cause its collapse, or that it can transform itself or reorganize itself. This is another formulation of Vergnaud's postulate that problems (perturbed system) are the sources and the criteria of knowing (Vergnaud 1981 p.220). It is important to realize that nothing is said about the process leading to the recovery of the equilibrium under the said constraints. They are proscriptive (Stewart, 1994 pp. 25-26), which means that they express necessary conditions to ensure the system viability, but not prescriptive, which means that they do not say in what way an equilibrium must be recovered.

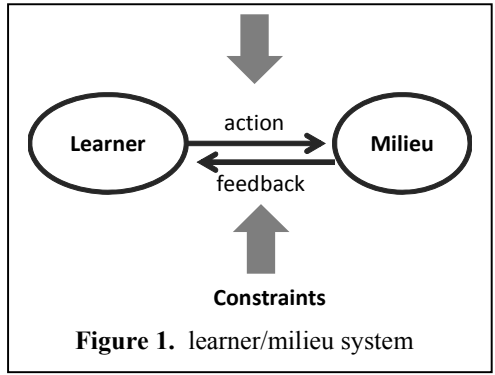

**Figure 1.** learner/milieu system

Hence, a definition of conception:

*A conception is the state of dynamical equilibrium of an action/feedback loop between a learner and a milieu under proscriptive constraints of viability.*

The study and characterization of a conception will be based on observable behaviors of the system (action, feedback) and outcomes of its functioning. It requires evidence of the assessment of the equilibrium, which depends on the possibility to elicit the learner's control of the interaction and of the milieu's reification of failures and success by adequate feedback.

Geometry provides many good examples: constructing a diagram on a sheet of paper with a pencil is permissive to empirical adjustments, while dynamic geometry software allowing messing up a diagram by dragging points can reify the failure of conforming to geometrical properties (Healy et al. 1994)—but still, "students may modify the figure 'to make it look right' rather than debug the construction process" (Jones1999 p.254).

The situated nature of a conception means that for different situations considered conceptually the same or for problems claimed isomorphic, one may associate different conceptions with the same learner. There is a large documentation of this phenomenon in the literature, for example under the theme of transfer, or from research in ethnomathematics. Anyhow, in the researcher's referential system, these different states of the observed systems [learner in a situation] should be labeled in the same way. For this reason, I define a *learner's knowing* as the set of conceptions which can be activated by different situations the observer considers conceptually the same—a qualification that indeed one will have to clarify. I realize that using "knowing" as a noun is rare, but it helps keeping distance with the word "knowledge" which has in education a strong authoritative connotation.

Having this definition of conception and of its relation to knowing, in the next section, I propose a model inspired by Vergnaud's formalization which I develop in the next section. Some examples will illustrate the model and facilitate clarifying the intention of my research program.

**Outlines of a model**

What is a model could be discussed at length within the PME multidisciplinary research community. I take here a pragmatic position, looking for something likely to facilitate our collaboration with research in educational technology, but also as a means to make more efficient our own research. The objective is to contribute to a better understanding of learners' understanding and to have eventually a practical value for teachers and designers.

It should be emphasized that the terms "conception", "knowing", "concept" and several others appearing in the description of the model are abstract terms, whose meaning is that of their functions and relations within the model. Indeed, we must then discuss how far the proposed formalization makes sense when confronted with "reality", and if it is an adequate tool for our research. The examples I will present and some additional comments will hopefully partly respond to this preoccupation.

I call "conception" a quadruplet $(P, R, L, \Sigma)$ in which:
- P is a set of problems.
- R is a set of operators.
- L is a representation system.
- $\Sigma$ is a control structure.

The first three elements are almost directly borrowed from the Vergnaud triplet. The vocabulary is different to avoid confusion with the vocabulary of psychology. In particular, operators correspond to actions one can observe in the functioning of the learner/milieu system; they are not schemes in psychological terms. The representation system is formed of all the semiotic tools which allow representing problems, supporting interaction and representing operators if formulations were required. The characterization of the set of problems P is more complex than expected. Two opposite solutions have been proposed: (i) to include all problems for which the conception provides efficient tools (Vergnaud 1991 p.145), but for basic concepts this option is too general to be effective; (ii) to consider a finite set of problems from which other problems will derive (Brousseau 1997 p.30), but this option opens the question of establishing that such a generative set of problems exists for any conception. A solution familiar to most researchers consists of deriving the description of P from both the observation of students in situations and from the analysis of historical and contemporary practices of mathematics. Actually, what one does when working on specific conceptions is to open a window on P by making explicit a few good representatives of its potential elements. These representatives work as kind of prototypical problems; this is a pragmatic implementation of Brousseau's proposal.

The forth element of the quadruplet, the control structure, includes behaviors such as making choices, assessing feedback, making decisions, judging the advancement of a problem solving process. These metacognitive behaviors are more often than not silent and invisible, hence rarely accessible to observation. It is why, to overcome this difficulty, one uses specific experimental settings, for example inviting learners to work in pairs, with the expectation that this will be enough to elicit these behaviors.

It is worth noticing that the quadruplet is not more related to the learner than to the milieu with which he or she interacts: the representation system allows the formulation and the use of the operators by the active sender (the learner) as well as the reactive receiver (the milieu); the control structure allows expressing the learner's means to assess an action, as well as the criteria of the milieu for selecting a feedback. It is in this sense that the quadruplet characterizing a conception is congruent to the conceptual definition of a conception as a property of the learner/milieu system.

This formalization not only allows characterizing conceptions and hence providing a framework to discuss their diagnosis, it has also the potential of helping to establish links among conceptions more precisely.

## Shaping relations between conceptions

### Arithmetic, from fingers to keystrokes

Addition has been widely studied, so there is enough resources to document what learners' conceptions could be like. This first example shows how eliciting the four dimensions of the quadruplet provide a synthetic and precise picture of the conceptions chosen for the purpose of the illustration (there are several others).

> **C1: Verbal counting IIIII & IIII**
> **P** -- *prototype*: "You have 5 pebbles, I give you 4 more, how many have you now?" (Objects are present or represented in an analogical way, both numbers are small).
> **R** – match fingers and objects, match fingers and number names, pointing to objects,
> **L** -- body language (finger counting, pointing), number naming, verbal counting
> $\Sigma$ – not counting twice an object, counting all the objects, order of the number names
> **C 2: Counting on 16 & 4**
> **P** -- *prototype*: "You have 16 pebbles, I give you 4 more, how many do have you now?" (The numbers are given, but the collections are not present, one of the numbers must be small enough)
> **R** -- choose the greater number, count on to determine the result.
> **L** -- body language (finger counting), number naming, verbal counting.
> $\Sigma$ -- order of the number names , match of fingers to number names
> **C3: written addition 16+23**
> **P** -- adding two integers
> **R** -- algorithm of column addition
> **L** -- decimal representation of numbers
> $\Sigma$ -- check the implementation of the algorithm, check the layout of number addition
> **C4: Pocket calculator [1][6][+][2][3][=]**
> **P** -- adding two integers, the result is bound by the size of the screen
> **R** – keystroke to represent a number, to process number addition
> **L** -- body language (keystrokes), decimal representation of numbers on the screen
> $\Sigma$ – keystrokes verification, order of magnitude.

**Table 1:** Four examples of conceptions of addition

First, the quadruplet puts to the fore the domain of validity of these conceptions, none of which can be claimed wrong but might be badly adapted outside their spheres of practice. For example, C1 will not be reliable in the domain of C3 and really difficult to implement in that of C2. For very large numbers, C4 will not work unless extended with additional strategies to deal with the limits imposed by the technology (e.g. screen display). Beyond the remark that the two first concern quantities and the two others concern numbers, one can compare the conceptions on each dimension of the quadruplet and analyze in an accurate way their commonalities and differences. One can also express relations often considered among conceptions; for this it is necessary to introduce a function allowing passing from one system of representation to another. Let us take the case of generality:

[*C is more general than C' if there exists a function of representation*

$f: L' \to L$ *so that* $\forall p \in P', f(p) \in P$]

So, one can show that C3 is more general than C2 and C4, but for different reasons although in both cases the size of the numbers is at stake. Indeed, we obtain in this way something that intuitively and informally would have been seen without such a sophisticated formulation. Isn't there here a flavor of pedantry? I hope not. It is just good to check that the formalism can express familiar facts but, indeed, it must go far beyond that and be capable of revealing less obvious relations.

**The challenge of translation**

The exercise of formalizing generality puts on the fore the importance of *translation*, a process we engage anytime we analyze learners activities. To emphasize the importance of the manipulation of representation systems, I take another example from the history of mathematics: ancient Egyptian arithmetic. The mathematical papyruses have been translated in the contemporary mathematical language, including the processes used to solve some problems. Figure 2 displays the sequence of steps to compute "10 times 1/5" (Couchoud 1993 pp.21-22). I will not here explain how the scribes obtained the result, 2, out of this sequence and how they moved from one line to the other. Looking closer to these processes suggests that the translation of ⟨symbol⟩ by 1/5 is misleading. What is denoted by the Egyptian sign is "five parts of the whole", hence an integer but integers which could not be added as integers are. These representations were not computed, instead scribes used tables to establish the correspondence between two numbers to be multiplied and the result. The control structures associated to this ancient conception of "fractions" and the modern one are completely different. Moreover, if one wants to consider passing back from the modern conception to the Egyptian one, that is expressing any fraction as a sum of unitary fractions, he or she will enter a new mathematical chapter. Several algorithms are available to compute an Egyptian decomposition for any fraction. For example, for 4055/4093 one will get the shortest and unique additive decomposition: [1/2 + 1/3 + 1/7 + 1/69 + 1/30650 + 1/10098761225]. Unfortunately, Egyptians could not write the last term. Analyzing these conceptions along the dimensions of the quadruplet in a systematic way makes easier figuring out what separates them.

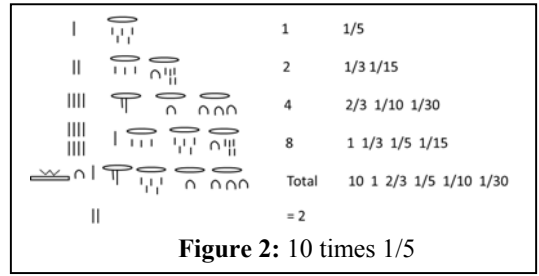

**Figure 2:** 10 times 1/5

Discussing learners' conceptions in the context of the mathematical curriculum is more difficult because operators and representation systems are often very similar. In this case, the control structure may be the discriminating element. As a matter of fact, this touches the foundation of conceptions because of the legitimacy controls provide by validating them.

Let us take a case in school algebra. In his research questioning the "production of meaning for Algebra", Romulo Lins (2001 p.47) observes the activity of students to whom he proposed the following task:

*To calculate how many oranges will fill into each box, we divide the total number of oranges by the number of boxes, i.e.:*

$$\text{orange per box} = \frac{\text{number of oranges}}{\text{number of boxes}}$$

*If I tell you the total number of oranges, and the number of oranges in each box, how would you calculate the number of boxes used?*

Justifying the task, Lins writes: "*The reason for presenting the 'algebraic' formula was to ascertain whether the pupils would constitute it into an object, dealing with it in the process of solving the problem; neither of them made any reference whatsoever to this formula*" (ibid.) In a very pragmatic way students manipulated oranges and boxes: "*They always used a number of something*" (ibid.). They dealt with quantities and not numbers. The control on their reasoning comes from the concrete reference the context makes possible. Actually, if algebra had been called up, literals would have been used to speak about actions on objects of a referent world made of boxes and oranges. This phenomenon is familiar, as Boero (2001 p.108) reports following a research he carried out in a different context: "*some students seem to transform the problem situation by thinking about the number of sheet and the weight of the envelope as physical variables*" whereas others "*put into a numerical equation the problem situation and transform the equation*".

Invited to write a postscript to the book "*Perspective on School Algebra*" (Sutherland et al. 2001), where I analyzed the reports of Lins, Boero and several others, I introduced the expression "*symbolic arithmetic*" to distinguish from algebra those conceptions in which symbols are manipulated and rules used with the supervision of a control structure grounded in the referent context.

In the following case, the role of the control structure clearly makes a difference between two conceptions which on the other hand seem to share the representation systems and operators.

**Questioning controls to understand representations**

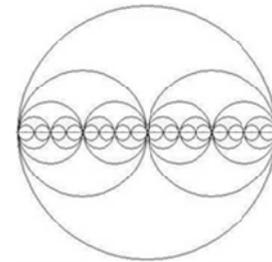

This case, which I have often presented, is as an excellent prototypical example of the complexity of identifying conceptions. The excerpt is borrowed from the work of Bettina Pedemonte (2002) on argumentation, cognitive continuity and proof. She chosen the problem described below and proposed it to pairs of students—the idea of having pairs was driven by the expectation to get spontaneous comments about actions, choices and decisions taken during the problem solving process.

> "Construct a circle with AB as a diameter. Split AB in two equal parts, AC and CB. Then construct the two circles of diameter AC and CB… and so on.
>> How does the perimeter vary at each stage?
>> How does the area vary?"

Using the formulas they know well both students, Vincent and Ludovic, express the perimeter and the area for the first steps in the series of drawings. They agree to conjecture that the perimeter will be constant and that the area will decrease to zero. But soon Vincent notices that "*the area is always divided by 2…so, at the limit? The limit is a line, the segment from which we started…*" This observation raises a conflict about the value of the perimeter which, in the opinion of Vincent, should be the length of the segment:

41. *Vincent*: It falls in the segment… the circles are so small.
42. *Ludovic*: Hmm… but it is always $2\pi r$.
43. *Vincent*: Yes, but when the area tends to 0 it will be almost equal…
44. *Ludovic*: No, I don't think so.
45. *Vincent*: If the area tends to 0, then the perimeter also… I don't know…
46. *Ludovic*: I will finish writing the proof.

Although Vincent and Ludovic collaborate well and seem to share the mathematics involved, the types of control they have on their problem-solving activity differ. Ludovic is working in the algebraic setting (Douady 1985) where control is based on a constant checking of the correctness of the symbolic manipulations conforming to the syntax of elementary algebra. Vincent is working in a symbolic-arithmetic setting where the control comes from a constant confrontation between what the formula "tells" and what is displayed by the drawings. So, both understand the initial situation in the "same" way, both manipulate the symbolic representations (i.e., the formulas of the perimeter and of the area) following the right syntax, but their controls are radically different. The symbolic representation supports the cooperation of the problem-solvers but it does not impose a shared understanding: as a *boundary object* it is flexible enough to adapt to the different meanings but robust enough to work as a tool for both students. To identify the differences beyond the apparent commonalities of representations one has to question learners' decisions and choices, which means identifying the control grounding their activity.

## Conception, knowing and concept

Understanding learners' conceptions requires their interpretation from the perspective of our own conception which we claim related to the same content of reference; one may say: the same concept. This can be expressed within the terms of the model, putting on the fore the role of translation which is more often than not implicit in our research practice. Let's take the case of "falsity" which is defined in the model in the following way (with a natural coding of the respective quadruplets):

> [*C is false from the point of view of C' if there exists a function of representation $f: L \rightarrow L'$,*
>
> *and there exists [$p \in P$, $r \in R$, $\sigma \in \Sigma$, $\sigma' \in \Sigma'$] so that $\sigma(r(p))=true$ and $\sigma'(f(r(p)))=false$*]

In other words, there exists a problem from the sphere of practice of C which has an accepted solution but which is assessed "false" from the point of view of C'.

"Generality" and "falsity" are not properties of conceptions but relations between two conceptions whose validity depends on the translation from one system of representation to the other. This is a general situation often hidden by the fact that we tend to read the production and the processes learners carry out directly in mathematical terms. Not being aware of this may make understanding learners difficult, as illustrated by Linns remarks on students' inability to escape the concrete reference of a situation he sees himself as mathematical. More generally, we have a tendency, often implicit, to consider ourselves as privileged knowers entitled to judge and evaluate other people's knowings. Such ambition requires at least that we can claim that the conception one assesses and the conception one holds are—so to say—ontologically compatible; they are concerned by the same object. This is difficult in mathematics where the only tangible things one manipulates are representations, and representations of representations. This can be solved within the model, taking Vergnaud's postulate as a grounding principle: problems are sources and criteria of knowings (1981 p.220):

> *Let C and C' be two conceptions and $C_a$ be a conception more general than C and C'. This means, with a natural coding of the respective quadruplets, that there exist functions of representation $f: L \rightarrow L_a$ and $f': L' \rightarrow L_a$ which relate C, C' and $C_a$. Then:*
>
> [*C and C' have the same object with respect to $C_a$ if for all p from P there exists p' from P' such that $f(p)=f'(p')$, and reciprocally*]

The fact that two conceptions have the same object does not mean that they have another type of relationship (one being false with respect to the other, or more general, or partial, or else), it may be the case that problems of P' (resp. P) cannot be expressed with L (resp. L'); and if they are, the translated problems may not be part of the sphere of practice of the other conception (e.g. the case of the conceptions of addition, Table 1). Eventually, conceptions have the same object if their defining problems (or their spheres of practice) can be matched from the point of view of a more general conception which in our case is the conception of the researcher. Research on mathematics learning must start with an introspection of researchers' own conception of the content at stake; questioning this conception is the first methodological tool to assess the validity of what can be said about learners' conception. This corresponds to the a priori analysis in the methodology associated to the TDS.

"To have the same object with respect to a conception $C_a$" sets an equivalence relation among conceptions. Let's now claim the existence of a conception $C\mu$ more general than any other conception to which it can be compared; this seems to be an abstract declaration, but pragmatically it corresponds to a piece of a mathematical theory.

*I call "concept" the set of all conceptions having the same object with respect to $C\mu$.*

This definition is aligned with the idea that a mathematical concept is not reduced to the text of its formal definition, but is the product of its history and of the set of practices in different communities. Indeed, there is no agent holding the concept and no way to ensure that we can enumerate a complete list of these conceptions. So, a last definition will allow reducing the distance between this abstract definition and the needs we have to have a practical model:

*I call "knowing" any subset of a concept which can be ascribed to a cognitive subject or a community.*

In practical terms, this definition of conception and knowing provides a framework which preserves learners' epistemic integrity despite contradictions and variability across situations. In a way which might seem more usual, I could summarize the ideas presented here in the following way: a *conception* is the instantiation of a *knowing* by a situation (it characterizes the subject/milieu system in a situation), or a *conception* is the instantiation of a *concept* by a pair (subject/situation).

The name cK¢ comes from the names of the three pillars of the model: *conception*, *knowing*, *concept*. I keep the word "knowledge" to name a conception which is identified and formalized by an institution (which is a body of an educational system in our case).

## Problems, conceptions and learning

Indeed, most problems are not solved by activating just one conception. So we need to be able to express the relation between a problem and conceptions which contribute to its solution, but without having to give details about this solution because of the cost of a too thin granularity. For this reason, I propose the more general idea of a *set of conceptions solving a problem*:

*Let p be a problem, and $\{C_1,…, C_n\}$ a set of conceptions.*

*$\{C_1,…, C_n\}$ solves p iff there exists a sequence of operators $(r_{i1}, …, r_{im})$ whose terms are taken in one of the $R_i$ so that the sequence $(p_1=r_{i1}(p), … , p_{im}=r_{im}(p_{im-1}))$ verifies that there exists σ from $\Sigma_{im}$ so that $\sigma(p_{im})$=solved.*

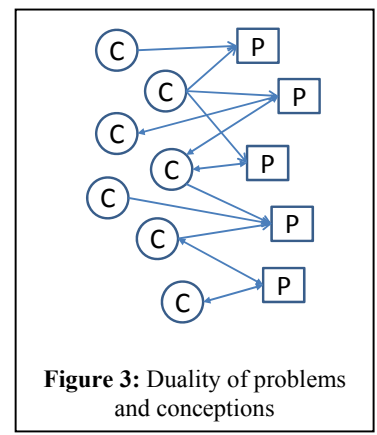

**Figure 3:** Duality of problems and conceptions

From this general characterization, one can derive more precise properties, for example expressing that a conception is specific to a problem (any set of conceptions solving it contains this conception), or that conceptions are equivalent from a problem solving perspective (one can replace the other without changing its property vis-à-vis the problem). Exploring this possibility evidences that problems and conceptions are of a dual nature: on the one hand conceptions need problems as constituents of their characterization, and on the other hand problems get their meaning from the conceptions contributing to their solutions. This duality suggests a natural connection between conceptions by the mediation of problems: this is exactly the idea of Vergnaud's conceptual field.

Learning is a process whose outcome is an evolution of conceptions being reinforced, questioned or transformed. The motor of this process is problems, which are (in our terms) destabilizations of the learner/milieu system. This destabilization can be obtained by modifications of constraints on the interaction between the learner and the milieu or modifications of the characteristics of the milieu (cf. Figure 1)—indeed, the learner is a "subsystem" on which no direct action is possible. The most difficult task is to find problems questioning the control structure and/or the representation system because the former is mostly implicit in the activity of learners and the role of the latter is invisible to their eyes once they are familiar with it. To overcome these difficulties is the *raison d'être* of the TDS situations of formulation and validation, which are based on a social organization of the class and a play on the characteristics of the milieu. I will not elaborate on this relation between cK¢ and the TDS, a part from noticing briefly here that within a didactical *problématique*, learning can be modeled as a transition function on a bipartite graph of conceptions and problems. Problems are the means to activate and (i) diagnose a conception, (ii) destabilize a conception to obtain an evolution, (iii) reinforce a conception.

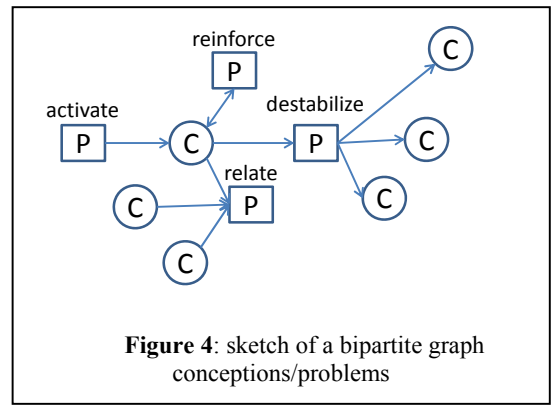

**Figure 4**: sketch of a bipartite graph conceptions/problems

Let C be the current diagnosed conception and $C_t$ a targeted conception (the expected learning outcome). The most critical evolution is the one to be obtained when C is false with respect to $C_t$. To engage the learning process, it is necessary to find a problem for which a representation is possible within both C and $C_t$, and which could be a means to reveal a conflict: a solution is conceivable from the perspective of C, but is not accessible in a way that the milieu witnesses and the subject recognizes. It may appear that such a problem does not exist and that intermediary problems, and possibly intermediary conceptions, are necessary to "reach" $C_t$. Learning is often not a single step but a path in the graph. To identify and create the conditions to bring this path to reality within a learning situation, in particular within a classroom, is one of the core objectives of the TDS.

Let's take the classical and well documented case of the sum of the angles of a triangle. The most common initial conception is that the bigger the object the bigger its measurable characteristics (e.g. area, perimeter), hence the sum of the angles. The operators are those involved in the manipulation of geometrical instruments (rule, compass, protractor, etc.) and symbolic arithmetic, the control structure includes visual control of actions and checking of computations. To activate this conception, one can ask students to draw triangles, measure angles and add up the obtained results. The variety of the results in the class will not raise questions since triangles are different; students will be reinforced in the confidence that they have the capacity to achieve the task. But measuring angles is not sufficient to destabilize the conception and give room to the targeted conception which is rooted not in the manipulation of the geometrical "object" but of its property. A possible way out is of asking learners to repeat the task, all with the same triangle (e.g. reproduced on a worksheet); the problem of deciding of the results for *this* single triangle will emerge and challenge the operators of the initial conception. More is needed to question its core theorem-in-action, this will come from the orchestration of the confrontation of the outcomes of collective workshops on triangles with contrasting shapes (small, large, sharp, flat). The destabilization of the initial conception can be overcome only by engaging in geometrical arguments of a theoretical nature. At this point, operators and controls are questioned, but the ambivalent nature of the triangle being both an object (of the spatio-graphic space) and a (geometrical) representation will (probably) remain unsolved.

## Concluding forewords

cK¢ proposes a framework for "learners modeling" taking up the challenge of providing a model of an epistemic relevance to bridge research in mathematics education and research on educational technology. It responds to a need for representations both understandable by researchers in mathematics education and computationally tractable. The formalism it dares should enhance the way one informs the design of technology enhanced learning environments, complementing descriptions generally available in natural language with no standardized narrative structure.

A European multidisciplinary assessment project (Baghera 2003) has been an occasion to check the efficiency of cK¢ in supporting a productive conversation between researchers in education and in computer-science. But probably more interesting for us is the powerfulness of this framework to think and develop our own research. Research in mathematics education develops jointly theories and experimentations, in this context models serve as mediators between theories of which they require an articulate and precise understanding, and experiments of which they frame the design and drive the collection of data. However, both theories and experiments raise difficult issues. On the side of theories, one has to deal with a complex discourse which rarely makes explicit all details and hence gives room to non-univocal interpretations. On the side of experiments, the practical implementation is always richer and more complex than what the design of models anticipates. Moreover, in the case of conceptions, one is confronted with issues (that Toulmin already noticed when proposing a model of argumentation): distinguishing operators from controls is not absolute (e.g. theorems can be activated as tools or predicates), controls are more often than not implicit. Such difficulties require further theoretical as well as methodological investigations.

Initially based on the Theory of Didactical Situation and the Theory of Conceptual Field, the cK¢ modeling framework is not restricted to them. For the purpose of its development and in order to enhance its efficiency it is necessary to integrate other theories to strengthen its components (e.g. representation, control system). But cK¢ holds other promises; it facilitates building a bridge between knowing and proving, constructing a link between control and proof, hence facilitating understanding the relation between argumentation and proof. But this is another topic which connects the research agendas I have had along my career, first on the teaching and learning of mathematical proof, then modeling learners' conceptions for the design of learning environments.

## Acknowledgements

I would like to express my thanks to Professor Patricio Herbst, Professor Anna Sierpinska and Dr. Bettina Pedemonte for their comments on an earlier draft of this paper, and Professor R. Knott for his "Calculator to convert a fraction into an Egyptian fraction" [http://www.maths.surrey.ac.uk/hosted-sites/R.Knott/Fractions/egyptian.html] retrieved 06-06-2013.

## Elements of a bibliography